\documentclass[11pt]{amsart}

\usepackage{epsfig, subfig, amsmath, amssymb, latexsym, amsfonts, xcolor}
\usepackage{mathdots}
\newcommand\dss{\displaystyle}

\usepackage{latexsym, amsfonts, amssymb}

\newenvironment{thm}{\subsection{}{\textbf {Theorem.}}\em}{}
\newenvironment{prop}{\subsection{}{\textbf {Proposition.}}\em}{}
\newenvironment{cor}{\subsection{}{\textbf {Corollary.}}\em}{}
\newenvironment{lem}{\subsection{}{\textbf {Lemma.}}\em}{}

\newenvironment{pf}{\noindent{\textbf {Proof.}}}
{\begin{flushright}\eop \end{flushright}\smallskip}

\newenvironment{defn}{\subsection{}{\textbf {Definition.}}\em}{\smallskip}

\newenvironment{ques}{\subsection{}{\textbf {Question.}}}{\smallskip}

\newcommand\fF{\ensuremath{\mathfrak F}}

\newcommand\fM{\ensuremath{\mathfrak M}}

\newcommand\fX{\ensuremath{\mathfrak X}}

\newcommand\cB{\ensuremath{\mathcal B}}
\newcommand\cC{\ensuremath{\mathcal C}}

\newcommand\cH{\ensuremath{\mathcal H}}

\newcommand\cK{\ensuremath{\mathcal K}}

\newcommand\cM{\ensuremath{\mathcal M}}

\newcommand\cP{\ensuremath{\mathcal P}}

\newcommand\cY{\ensuremath{\mathcal Y}}

\newcommand\bbC{\ensuremath{\mathbb C}}

\newcommand\bbM{\ensuremath{\mathbb M}}
\newcommand\bbN{\ensuremath{\mathbb N}}

\newcommand\bbT{\ensuremath{\mathbb T}}

\newcommand\bbZ{\ensuremath{\mathbb Z}}

\newcommand\hilb{\ensuremath{\mathcal H}}

\newcommand\ol{\ensuremath{\overline}}

\newcommand\eop{{{\hfil \ensuremath \Box}}}

\newcommand\norm{\ensuremath {\Vert}}

\newcommand\bofh{\ensuremath{\cB ( \cH)}}

\begin{document}
\title{Normal operators with highly incompatible off-diagonal corners}

\thanks{${}^1$ Research supported in part by NSERC (Canada)}
\thanks{${}^2$ Research was supported in part by  Natural Science Foundation
for Young Scientists of Jilin Province (No.: 20190103028JH),
National Natural Science Foundation of China (No.: 11601104, 11671167, 11201171) and China Scholarship Council (No.:201806175122).}

\thanks{{\ifcase\month\or Jan.\or Feb.\or March\or April\or May\or
June\or
July\or Aug.\or Sept.\or Oct.\or Nov.\or Dec.\fi\space \number\day,
\number\year}}
\author
	[L.W.~Marcoux]{{Laurent W.~Marcoux${}^1$}}
\address
	{Department of Pure Mathematics\\
	University of Waterloo\\
	Waterloo, Ontario \\
	CANADA  \ \ \ N2L 3G1}
\email{LWMarcoux@uwaterloo.ca}

\author
	[H.~Radjavi]{{Heydar Radjavi}}
\address
	{Department of Pure Mathematics\\
	University of Waterloo\\
	Waterloo, Ontario \\
	CANADA  \ \ \ N2L 3G1}
\email{HRadjavi@uwaterloo.ca}

\author
	[Y.H.~Zhang]{{Yuanhang~Zhang${}^2$}}
\address
	{School of Mathematics\\
	Jilin University\\
	Changchun 130012\\P.R.CHINA}
\email{zhangyuanhang@jlu.edu.cn}

\begin{abstract}
Let $\hilb$ be a complex, separable Hilbert space, and $\bofh$ denote the set of all bounded linear operators on $\hilb$.   Given an orthogonal projection $P \in \bofh$ and an operator $D \in \bofh$, we may write $D=\begin{bmatrix} D_1& D_2 \\ D_3 & D_4 \end{bmatrix}$ relative to the decomposition $\hilb = \mathrm{ran}\, P \oplus \mathrm{ran}\, (I-P)$.    In this paper we study the question:  for which non-negative integers $j, k$ can we find a normal operator $D$ and an orthogonal projection $P$ such that $\mathrm{rank}\, D_2 = j$ and $\mathrm{rank}\, D_3 = k$?  Complete results are obtained in the case where $\mathrm{dim}\, \hilb < \infty$, and partial results are obtained in the infinite-dimensional setting.
\end{abstract}


\maketitle
\markboth{\textsc{  }}{\textsc{}}



\section{Introduction}

\subsection{} \label{sec1.01}
Let $\hilb$ denote a complex, separable Hilbert space.  By $\bofh$ we denote the space of bounded linear operators acting on $\hilb$, keeping in mind that when $\mathrm{dim}\,\hilb = n < \infty$ we may identify $\hilb$ with $\bbC^n$, and $\bofh$ with $\bbM_n(\bbC)$.    We write $\cP(\hilb) := \{ P \in \bofh: P = P^2 = P^*\}$ to denote the set of orthogonal projections in $\bofh$.    Given $T \in \bofh$, $T$ admits a natural $2 \times 2$ operator-matrix decomposition 
\[
T = \begin{bmatrix} T_1 & T_2 \\ T_3 & T_4 \end{bmatrix} \]
with respect to the decomposition $\hilb = P \hilb \oplus (I-P) \hilb$.   
Of course, $T_j = T_j(P)$, $1 \le j \le 4$.  

We are interested in determining to what extent the set $\{ (T_2(P), T_3(P)) : P \in \cP(\hilb)\}$ determines the structure of the operator $T$.   Following~\cite{LMMR2018}, we say that $T$ has property (CR) (the \textbf{common rank} property) if $\mathrm{rank}\, T_2(P) = \mathrm{rank}\, T_3(P)$ for all $P \in \cP(\hilb)$.  We recall that an operator $A \in \bofh$ is said to be \textbf{orthogonally reductive} if $P \in \cP(\hilb)$ and $(I-P) A P = 0$ implies that $P A (I-P) = 0$.   That is, every invariant subspace for $A$ is orthogonally reducing for $A$.      In the above-cited paper, the following result was obtained:


\begin{thm} \label{thm1.02}~\cite[Theorem~5.8]{LMMR2018}
Let $\hilb$ be a complex Hilbert space and $T \in \bofh$.  If $T$ has property (CR), then there exist $\lambda, \mu \in \bbC$ and $A \in \bofh$ with $A$ either selfadjoint or an orthogonally reductive unitary operator such that $T = \lambda A + \mu I$.
\end{thm} 


\subsection{} \label{sec1.03}
In fact, if $\mathrm{dim}\, \hilb < \infty$, then the converse is also true~(\cite[Theorem~3.15]{LMMR2018}).   We note that every normal operator (and hence every unitary operator) acting on a finite-dimensional Hilbert space is automatically orthogonally reductive.   In particular, every operator $T$ that has property (CR) must be normal with spectrum lying either on a line or a circle, and when $\hilb$ is finite-dimensional, every such normal operator has property (CR).

\smallskip

Property (CR) was termed a ``\emph{compatibility}" condition on the off-diagonal corners of the operator $T$.   In this paper, we examine to what extent the off-diagonal corners of a normal operator $D$ may be ``\emph{incompatible}" in the sense of rank.   That is, writing $D = \begin{bmatrix} D_1  & D_2  \\ D_3 & D_4  \end{bmatrix}$ relative to $\hilb = P\hilb \oplus (I-P) \hilb$, we consider how large  
\[
| \mathrm{rank}\, D_2  - \mathrm{rank}\, D_3| \]
can get.   

More generally, our main result (Theorem~\ref{thm2.05} below) shows that if $\mathrm{dim}\, \hilb= n < \infty$ and $1 \le j, k \le \lfloor \frac{n}{2} \rfloor$, then there exist a normal operator $D$ and a projection $P$ such that $\mathrm{rank}\, D_2(P) = j$ while $\mathrm{rank}\, D_3(P) = k$.   If $\mathrm{dim}\, \hilb = \infty$ and if $0 \le j, k \le \infty$, then the same conclusion holds (Theorem~\ref{thm3.02}).

\smallskip

The infinite-dimensional setting also allows for certain subtleties which cannot occur in the finite-dimensional setting.  For example, if $\mathrm{dim}\, \hilb = n < \infty$, $D = \begin{bmatrix} D_1 & D_2 \\ D_3 & D_4 \end{bmatrix} \in \bofh$ is normal and $D_3 = 0$, then $D_2 = 0$.   Indeed, this is just a restatement of the fact that every normal matrix is orthogonally reductive.  This follows by observing that the normality of $D$ implies that 
\[
D_1^* D_1  - D_1 D_1^* = D_2 D_2^* - D_3^* D_3.   \]
Thus  $\mathrm{tr}(D_2 D_2^*) = \mathrm{tr}(D_3^* D_3)$, or equivalently $\norm D_2 \norm_2 = \norm D_3 \norm_2$, where $\norm \cdot \norm_2$ refers to the Frobenius (or Hilbert-Schmidt) norm.   
From this, $D_3= 0$ clearly implies that $D_2 = 0$.    We shall show that if $\hilb$ is infinite-dimensional, then it is possible to have $D_3(P) = 0$ while $D_2(P)$ is a quasiaffinity (i.e. $D_2(P)$ has trivial kernel and dense range), although it is not possible for $D_3(P)$ to be compact and $D_2(P)$ to be invertible (see Proposition~\ref{prop3.03} below).


\subsection{} \label{sec1.04}
It is worth mentioning that a related question where ranks are replaced by unitarily invariant norms has been considered by Bhatia and Choi~\cite{BC2006}.   More specifically, they consider normal matrices $D = \begin{bmatrix} D_1 & D_2 \\ D_3 & D_4 \end{bmatrix}$ acting on $\hilb := \bbC^n \oplus \bbC^n$.  As noted above, normality of $D$ shows that $\norm D_2 \norm_2 = \norm D_3 \norm_2$.    In the case of the operator norm $\norm \cdot \norm$, it follows that $\norm D_3 \norm \le \sqrt{n} \, \norm D_2 \norm$, and equality can be obtained in this expression if and only if $n \le 3$.    (If we denote by $\alpha_n$ the minimum number such that $\norm D_3 \norm \le \alpha_n \, \norm D_2 \norm$ for all $D \in \bbM_{2n}(\bbC)$ as above -- so that $\alpha_n \le \sqrt{n}$ --  it is not even known at this time whether or not the sequence $(\alpha_n)_n$ is bounded.

It is interesting to note that the example they give for the case where $n = 3$ and $\alpha_3 = \sqrt{3}$ is also an example of a normal matrix $D \in \bbM_6(\bbC)$ for which $\mathrm{rank}\, D_2 = 1$ and $\mathrm{rank}\, D_3 = 3$.   


\vskip 1 cm

\section{The finite-dimensional setting}


\subsection{} \label{sec2.01}
In examining the incompatibility of the off-diagonal corners of a normal operator $D \in \bbM_n(\bbC)$, we first dispense with the trivial cases where $n \in \{ 2, 3\}$.   Indeed, as seen in Proposition~3.7 of~\cite{LMMR2018}, in this setting, $D$ automatically has property (CR). 

For this reason, henceforth we shall assume that $\mathrm{dim}\, \hilb \ge 4$.  

The key to obtaining the main theorem of this section is Theorem~\ref{thm2.03}, which shows that if $\mathrm{dim}\, \hilb = 2 m$ for some integer $m \ge 2$, then we can find a normal operator $D$ such that $\mathrm{rank}\, D_3 = 1$ and $\mathrm{rank}\, D_2 = m$.   For $m = 2$, this is an immediate consequence of Theorem~3.15 of~\cite{LMMR2018}, since in this case, given a normal operator $D \in \bbM_4(\bbC)$ whose eigenvalues do not lie either on a common circle or a common line, $D$ fails to have property (CR), and this can only happen if there exists a projection $P \in \bbM_4(\bbC)$ of rank two such that $\mathrm{rank}\, D_2(P) = 2$, while $\mathrm{rank}\, D_3(P) = 1$.  

Given $X = [x_{i,j}], Y=[y_{i,j}] \in \bbM_n(\bbC)$, we shall denote by $X \bullet Y$ the \textbf{Hadamard} or \textbf{Schur} product of $X$ and $Y$; i.e. $X \bullet Y = [x_{i,j}\,  y_{i,j}] \in \bbM_n(\bbC)$.


\begin{lem} \label{lem2.02}
Let $m \ge 3$ be an integer.   Let 
\[
A = \mathrm{diag} (\alpha_1, \alpha_2, \ldots, \alpha_m) \mbox{\ \ \ \ \ and \ \ \ \ } B = \mathrm{diag} (\beta_1, \beta_2, \ldots, \beta_m) \]
be diagonal operators in $\bbM_m(\bbC)$, and $D := \begin{bmatrix} A & 0 \\ 0 & B \end{bmatrix}$.   Set $Z := [ z_{j,k}] \in \bbM_m(\bbC)$, where $z_{j,k} := \alpha_j - \beta_k$ for all $1 \le j, k \le m$.  Suppose that there exists a positive definite matrix $S \in \bbM_m(\bbC)$ such that 
\[
\mathrm{rank}\, S \bullet Z = 1 \mbox{\ \ \ \ \ \ and \ \ \ \ \ } \mathrm{rank}\, S^t \bullet Z = m, \]
where $S^t$ denotes the transpose of $S$.   

Then there exists a projection $P \in \bbM_{2m}(\bbC)$ such that if $D = \begin{bmatrix} D_1 & D_2 \\ D_3 & D_4 \end{bmatrix}$ relative to $\bbC^{2m} = \mathrm{ran}\, P \oplus \mathrm{ran}\, (I-P)$, then $\mathrm{rank}\, D_2 = m$ and $\mathrm{rank}\, D_3 = 1$.
\end{lem}


\begin{pf}
We leave it as an exercise for the reader to show that  $0 < S \in \bbM_m(\bbC)$ implies that $S$ can be expressed in the form $S = M N^{-1}$, where $M$ and $N$ are  two commuting positive definite matrices satisfying $M^2 + N^2 = I_m$.  From this it follows that
\[
P:=
	\begin{bmatrix}
	M^2 & M N \\
	M N & N^2
	\end{bmatrix} \]
is an orthogonal projection in $\bbM_{2m}(\bbC)$ whose rank is $m = \mathrm{tr}(P)$.    Since $P = \begin{bmatrix} M \\ N \end{bmatrix} \ \begin{bmatrix} M & N \end{bmatrix}$, we deduce that $\begin{bmatrix} M \\ N \end{bmatrix}$ is an isometry from $\bbC^m$ into $\bbC^{2 m}$.  A straightforward computation shows that
\[
I_{2 m} - P = \begin{bmatrix} I_m - M^2 & - M N \\ - MN & I_m - N^2 \end{bmatrix} = \begin{bmatrix} N \\ - M \end{bmatrix} \ \begin{bmatrix} N & - M \end{bmatrix}, \]
and that $\begin{bmatrix} N \\ - M \end{bmatrix}$ is once again an isometry of $\bbC^m$ into $\bbC^{2 m}$.

\smallskip

Our goal is to show that $\mathrm{rank}\, (I-P) D P = 1$, while $\mathrm{rank}\, P D (I-P) = m$.   In light of the fact that both $\begin{bmatrix} M \\ N \end{bmatrix}$ and $\begin{bmatrix} N \\ - M \end{bmatrix}$ are isometries, this is equivalent to proving that
\[
\mathrm{rank}\, (N A M - M B N) = \mathrm{rank} \begin{bmatrix} N & -M \end{bmatrix} \ \begin{bmatrix} A & 0 \\  0 & B \end{bmatrix} \ \begin{bmatrix} M \\ N \end{bmatrix} = 1, \]
while
\[
\mathrm{rank}\, (M A N - N B M) = \mathrm{rank} \begin{bmatrix} M & N \end{bmatrix} \ \begin{bmatrix} A & 0 \\  0 & B \end{bmatrix} \ \begin{bmatrix} N \\ -M \end{bmatrix} = m. \]

Now $N$ and $M$ are each invertible in $\bbM_m(\bbC)$, and $N M = M N$ implies that $N^{-1}$ and $M$ also commute.  Thus
\[
\mathrm{rank}\, (N A M - M B N) = \mathrm{rank}\, (A M N^{-1} - N^{-1} M B) = \mathrm{rank}\, ( A S - S B) = \mathrm{rank}\, S \bullet Z = 1, \]
while
\begin{align*}
\mathrm{rank}\, (M A N - N B M) 
	&= \mathrm{rank}\, (N^{-1} M A - B M N^{-1}) \\
	&= \mathrm{rank}\, (S A - B S) \\
	&= \mathrm{rank}\, (A S^t - S^t B) \\
	&= \mathrm{rank}\, S^t \bullet Z\\
	&= m.
\end{align*}	
\end{pf}
	

\begin{thm} \label{thm2.03}
Let $m \ge 1$ be an integer.   Then there exist a normal operator $D \in \cB(\bbC^{2m}) \simeq \bbM_{2m}(\bbC)$ and an orthogonal projection $P$ of rank $m$ such that if $D = \begin{bmatrix} D_1 & D_2 \\ D_3 & D_4\end{bmatrix}$ relative to $\bbC^{2m} = \mathrm{ran}\, P \oplus \mathrm{ran}\, (I-P)$, then $\mathrm{rank}\, D_2 = m$ and $\mathrm{rank}\, D_3 = 1$.
\end{thm}


\begin{pf}
The case $m = 1$ is easily handled by the operator $D = \begin{bmatrix} 1 & 1 \\ 1 & 1 \end{bmatrix}$ and the projection $P = \begin{bmatrix} 1 & 0 \\ 0 & 1 \end{bmatrix}$.   The case where $m = 2$ follows from Proposition~3.13 of~\cite{LMMR2018}.

Suppose, therefore, that $m \ge 3$.  By Lemma~\ref{lem2.02}, we have reduced our problem to that of finding two diagonal matrices $A = \mathrm{diag} (\alpha_1, \alpha_2, \ldots, \alpha_m)$ and $B= \mathrm{diag} (\beta_1, \beta_2, \ldots, \beta_m)$,  and a positive definite matrix $0 < S = [s_{j, k}] \in \bbM_m(\bbC)$ such that
\[
\mathrm{rank}\, S \bullet Z = 1 \mbox{\ \ \ \ \ \ and \ \ \ \ \ } \mathrm{rank}\, S^t \bullet Z = m. \]

We begin by specifying $A$ and $B$;  we first temporarily fix a parameter $1 < \gamma$ whose exact value we shall determine later.  For $1 \le j \le m$, set $\alpha_j = j \gamma + i$.   Set $B=A^*$, so that $\beta_k = \ol{\alpha_k} = k \gamma - i$.   Then $Z = [z_{j, k}] = [ (j-k) \gamma + 2 i]$.

Next, we set $S (= S(\gamma)) =[s_{j, k}]$, where $s_{j, k} = \dss \frac{ 2i }{{ (j-k) \gamma + 2 i}}$.  Observe first that for $1 \le j, k \le m$,
\[
\ol{s_{k, j}} = \frac{-2i}{(k-j) \gamma -2i} = \frac{ 2i }{{ (j-k) \gamma + 2 i}} =  s_{j, k}, \]
so that $S$ is clearly hermitian, and $s_{j, j} = 1$ for all $1 \le j \le m$.  It is therefore reasonably straightforward to see that since $m$ is a fixed constant, and since $\dss \lim_{\gamma \to \infty} \dss \frac{2i}{(j-k) \gamma + 2i} = 0$ for all $1 \le j \ne k \le m$, there exists a constant $\Gamma (m) \ge 1$ such that $\gamma > \Gamma(m)$ ensures that $\norm S - I_m \norm < \frac{1}{4}$, and thus $S (=S(\gamma))$ must be positive definite.  

For an explicit estimate for $\Gamma(m)$, we may observe that if $R = [r_{j, k}] \in \bbM_m(\bbC)$, and if $\norm R \norm_\infty := \max \{ |r_{j,k}| : 1 \le j, k \le m \}$, then $\norm R \norm \le m \norm R \norm_\infty$.   Indeed, if $x=(x_k)_{k=1}^m \in \bbC^m$, then (using the Cauchy-Schwartz inequality) we find that 
\begin{align*}
\norm R x \norm^2 
	&= \sum_{j=1}^m | \sum_{k=1}^m r_{j,k} x_k |^2  \\
	&\le \sum_{j=1}^m  m \norm R \norm_\infty^2 \norm x \norm^2 \\
	&= m^2 \, \norm R \norm_\infty^2 \, \norm x \norm^2,
\end{align*}	
from which the result follows.
In particular, by choosing $\Gamma(m) = 8 m$, $\gamma > \Gamma(m)$ implies that
\[
\norm S - I_m \norm\leq m \max_{1\leq j,k\leq m}|s_{j,k}-\delta_{j,k}|=m \max_{1\leq j\neq k\leq m}|s_{j,k}| <m \frac{2}{\gamma}< \frac{1}{4}, \]
and so $S$ is a positive invertible operator.

Consider
\[
S \bullet Z = [s_{j, k}\  z_{j, k}] = [\frac{2i}{(j-k) \gamma + 2i} \ ((j-k)\gamma+2i)] = [2i]_{m \times m}. \]

It is clear that $S \bullet Z \in \bbM_m(\bbC)$ is a rank-one operator;  indeed, $S \bullet Z = 2m i Q$, where $Q$ is the rank-one projection whose matrix consists entirely of the entries $\frac{1}{m}$.

We therefore turn our attention to
\[
S^t \bullet Z = [s_{k, j} \ z_{j, k}] = [ \frac{2i}{(k-j) \gamma+2i} \ ((j-k)\gamma+2i) ] = [ 2i \theta_{j, k}], \]
where $\theta_{j, k} = \dss \frac{(j-k)\gamma+2i}{(k-j)\gamma + 2i} \in \bbT$, $1 \le j, k \le m$.  
Observe that if $1 \le j, k \le m-1$, then $\theta_{j, k} = \theta_{j+1, k+1}$.  Thus $T := \frac{1}{2i} (S^t \bullet Z)$ is  a Toeplitz matrix, and the diagonal entries of $T$ are all equal to $1$.

In fact, for $1 \le j, k \le m$,
\[
\ol{\theta_{k,j}} = \frac{(k-j) \gamma - 2i}{(j-k) \gamma - 2i} = \frac{- ((j-k)\gamma + 2i)}{-((k-j) \gamma + 2i)} = \theta_{j,k},  \]
and therefore $T$ is not only Toeplitz, but hermitian as well.

There remains only to show that the rank of $S^t \bullet Z$ is $m$, or equivalently, that $\mathrm{det}\, T \ne 0$.

\smallskip

Define $\hat{T}=2I_m-mQ$.
Then $\hat{T}$ is invertible and $\hat{T}^{-1}=\frac{1}{2-m}Q+\frac{1}{2}(I_m-Q)$.  Note that each diagonal entry of $\hat{T}$ is 1, 
while each off-diagonal entry is $-1$.
From this and the calculations above it follows that 
\[
\norm T - \hat{T} \norm\leq m\|T-\hat{T}\|_\infty=m(\underset{1\leq j\neq k\leq m}{\max}|\theta_{j, k}+1|) <m\frac{4}{\gamma}< \frac{1}{2}<\frac{1}{\|\hat{T}^{-1}\|}, \]
implying that $T$ is invertible, whenever $\gamma > \Gamma(m) = 8m$.

Thus, by choosing $\gamma >\Gamma(m) = 8m$, we see that a positive solution to our problem can be found.
\end{pf}


\subsection{} \label{sec2.04}
Suppose now that $n \ge 5$ is an integer and that $T \in \bbM_n(\bbC)$.   If $P \in \cP(\bbC^n)$ is any projection, then the minimum of $\mathrm{rank}\, P$ and $\mathrm{rank}\, (I-P)$ is at most $\lfloor \frac{n}{2} \rfloor$.  It follows that 
\[
\max(\mathrm{rank}\, T_2(P), \mathrm{rank}\, T_3(P)) \le \lfloor \frac{n}{2} \rfloor.\]
As already observed, if $D \in \bbM_n(\bbC)$ is normal, then $D$ is orthogonally reductive, and so if $\mathrm{rank}\, T_3(P) = 0$, then automatically $\mathrm{rank}\, T_2(P) = 0$.   In light of these observations, we see that the following result is the best possible, and it is the main theorem of this section.


\begin{thm} \label{thm2.05}
Let $n\geq 2$ be a positive integer, $1\leq j,k\leq \lfloor\frac{n}{2}\rfloor$.  Then there exist a normal operator $D \in \bbM_n(\bbC)$ and a projection $P$ such that relative to $\bbC^n = \mathrm{ran}\, P \oplus \mathrm{ran}\, (I-P)$, we can write
	\[
	D=\begin{bmatrix} D_1 & D_2 \\ D_3 & D_4 \end{bmatrix} \]
where $\mathrm{rank}\, D_2 = k$ and $\mathrm{rank}\, D_3 = j$.
\end{thm}


\begin{pf}
Without loss of generality, we can assume that $k\geq j$. First, we set $m := (k-j) +1$.   Applying Theorem~\ref{thm2.03} we may choose a normal element $M \in \bbM_{2 m}(\bbC)$ such that
\[
	M =\begin{bmatrix} M_1 & M_2 \\ M_3 & M_4 \end{bmatrix}, \]
where $\mathrm{rank}\, M_2 =(k-j)+1 $ and $\mathrm{rank}\, M_3 = 1$.
Define
\[
	\hat{D}=\begin{bmatrix} I_{j-1}&&&I_{j-1}\\ &M_1 & M_2 & \\ & M_3 & M_4 &\\I_{j-1}&&&I_{j-1} \end{bmatrix}. \]
Here, it is understood that if $j=1$, then $I_0$ acts on a space of dimension zero.   	
Finally, let 
\[
	D=0_{n-2k} \oplus \hat{D} = \begin{bmatrix} 0_{n-2k}&&&&\\ & I_{j-1} & & & I_{j-1}\\ & & M_1 & M_2 & \\ & & M_3 & M_4 & \\ & I_{j-1} & & & I_{j-1} \end{bmatrix}.\]
(Again, if $n = 2k$, the  $0_0$ term is not required.)	
Set $P = I_{(n-2k)+(j-1) + m} \oplus 0_{m+(j-1)}$, and relabel $D = \begin{bmatrix}D_1 & D_2 \\ D_3 & D_4\end{bmatrix}$ relative to the decomposition $\bbC^n = \mathrm{ran}\, P \oplus \mathrm{ran}\, (I-P)$.   It is then routine to verify that  $\mathrm{rank}\, D_2 = k$ and $\mathrm{rank}\, D_3 = j$.
\end{pf}


\subsection{} \label{sec2.06}
The operator $D$ constructed in Theorem~\ref{thm2.05} is far from unique.  Indeed, we first note that we were free to choose arbitrarily large $\gamma$'s in the definition of $A$ and $B$ defined above.     Secondly, it is not hard to show that by choosing $B= A^*$ and $Z$ as we did above, and by defining $S$ such that $S \bullet Z = 2i Q$, $S$ is always hermitian.   Thus, given one triple $(A, B, S)$ as above that works, if we slightly perturb the weights $\alpha_j$ of our given $A$ to obtain a diagonal matrix $A_0$ and we set $B_0 = A_0^*$, then the new $S_0$ we require to make $S_0 \bullet Z_0 = 2i Q$ will be sufficiently close to the original $S$ so as to be invertible (since the set of invertible operators is open in $\bbM_m(\bbC)$).   


\subsection{} \label{sec2.07}
An interesting, but apparently far more complicated question, is to characterise those normal operators $D \in \bbM_{2m}(\bbC)$ for which it is possible to find a projection $P$ of rank equal to $m$ such that $\mathrm{rank}\, (I-P) D P = 1$ and $\mathrm{rank}\, P D (I-P) = m$.    We are not able to resolve this question at this time.   We can assert, however, that not only is such a normal operator abstractly ``far away" from operators with property (CR);  in fact, we are able to quantify this distance, and say a bit more about the structure of $D$.

Let $n \ge 1$ be an integer, and recall that the function 
\[
\begin{array}{rccc}
\rho: & \bbM_n(\bbC) \times \bbM_n(\bbC) & \to & \{ 0, 1, 2, \ldots \} \\
	& (A, B) & \mapsto & \mathrm{rank}\, (A-B)
\end{array} \]
defines a metric on $\bbM_n(\bbC)$.   

We also recall that an operator $T \in \bofh$ (where $\mathrm{dim}\, \hilb \in \bbN \cup \{\infty\}$) is said to be \textbf{cyclic} if there exists $x \in \hilb$ such that $\mathrm{span}\, \{ x, Tx, T^2 x, \ldots \}$ is dense in $\hilb$.   Obviously this can only happen if $\hilb$ is separable, and it is well-known that a normal operator is cyclic if and only if it has \emph{multiplicity one}; that is, its commutant $N^\prime := \{ X \in \bofh: X N = N X \}$ is a \textbf{masa} (i.e. a maximal abelian selfadjoint subalgebra of $\bofh$).  If $N$ is a compact, normal operator, then this is equivalent to saying that the eigenspaces corresponding to the eigenvalues of $N$ are all one-dimensional, and together they densely span the Hilbert space.


\begin{thm} \label{thm2.08}
Let $m \ge 3$ be an integer, and suppose that $D \in \bbM_{2m}(\bbC)$ is a normal operator.   Suppose that $P \in \bbM_{2m}(\bbC)$ is an orthogonal projection of rank $m$ and that $D \in \bbM_{2m}(\bbC)$ is a normal operator for which $\mathrm{rank}\, (I-P) D P = 1$ and $\mathrm{rank}\, P D (I-P) = m$.

Then
\begin{enumerate}
	\item[(a)]
 	$D$ has $2m$ distinct eigenvalues (and therefore $D$ is a cyclic operator); and
	\item[(b)] 
	$\rho(D, Y)\geq \lfloor\frac{m-1}{2}\rfloor$ for all $Y \in \cY$, where $\cY$ is the set of matrices in $M_{2m}(\mathbb{C})$ which satisfy property (CR).
\end{enumerate}
\end{thm}


\begin{pf}
First observe that we may assume without loss of generality that $D$ is invertible, since otherwise we simply add a sufficiently large multiple of the identity to $D$, which  affects neither the hypotheses nor the conclusion of the Theorem.    
\begin{enumerate}
	\item[(a)]
	Next, we set $P_0 := P$, and let $V_0$ be the range of $P_0$.   By hypothesis,
	\[\mathrm{dim} (V_0\vee DV_0)=m+1,~~\mathrm{dim} (V_0\cap D^{-1}V_0)=m-1.\]

	More generally, we claim that the following chain of subspaces has strictly increasing dimensions (from $0$ to $n=2m$):
	 \[V_{-m}\subset V_{-m+1}\subset \cdots\subset V_{-1}\subset V_0\subset V_1\subset \cdots \subset V_m,\]
	 where
	 \[V_{k+1}=V_k\vee DV_k,~~\forall~~ 0\leq k\leq m-1,\]
	 and
	 \[V_{k-1}=V_k\cap  D^{-1}V_k,~~\forall~~ -m+1\leq k\leq 0.\]

	Assume to the contrary that this fails.   Let $P_k$ be the projection to the range of $V_k$, $-m\leq k\leq m$.
	\begin{enumerate}
		\item[(i)] 
		If $V_{k+1}=V_k$ for some $0<k<m$, then $DV_k=V_k$.  This implies that $D^*V_k=V_k$, i.e., that 
		$P_kD(I-P_k)=0$. Since $P_k\geq P_0$, we deduce that
		\[P_0D(I-P_k)=0,\]
		($\mathrm{rank}(P_k)=\mathrm{dim}(V_k)\leq m+k<2m$).
		In other words, $P_0D(I-P_0)$ has nontrivial kernel in $V_0$, a contradiction.
		\item[(ii)]
		Similarly, if $V_{k+1}=V_k$ for some $-m\leq k<0$, then once again $DV_k=V_k$ and $P_kD(I-P_k)=0$. Since $P_k\leq P_0$,
		we deduce that
		\[P_kD(I-P_0)=0.\]
		($\mathrm{rank}(P_k)=\mathrm{dim}(V_{k+1}) = \mathrm{dim}(V_k) \geq 1$.)
		This implies that the range $P_0D(I-P_0)$ is smaller that that of $P_0$; a contradiction.
	\end{enumerate}
	Thus the claim is proved.

	In particular, $V_{-m+1}$ is one-dimensional. Pick a unit vector in $V_{-m+1}$.
	We next show that $x$ is a cyclic vector for $D$.

	Note that $D x\notin V_{-m+1}$, and hence $x, Dx$ span $V_{-m+2}$.
	Under the assumption that $\{x,Dx,\cdots, D^jx\}$ spans $V_{-m+j+1}$, we see that  $\{x,Dx,\cdots, D^{j+1}x\}$ spans $V_{-m+j+2}$
	by construction. This is true for all $0\leq j\leq 2m-1$, which proves that $x$ is a cyclic vector of $D$.
	\item[(b)]
	With the decomposition of $\mathbb{C}^{2m}=\textup{ran}P\oplus \textup{ran} (I-P)$, we may write
	\[D=
		\begin{bmatrix}
		D_1 & D_2 \\
		D_3 & D_4
		\end{bmatrix}.\]
Next, suppose that $Y \in \cY$, so that $Y$ has the common rank property.
With respect to the same decomposition of $\mathbb{C}^{2m}$, we have that 
\[Y=
	\begin{bmatrix}
	Y_1 & Y_2 \\
	Y_3 & Y_4
	\end{bmatrix}. \]
Define $F := D - Y$ and write 
\[F =
	\begin{bmatrix}
	F_1 & F_2 \\
	F_3 & F_4
	\end{bmatrix}. \]
Clearly $D_2 = Y_2 + F_2$ and $D_3 = Y_3 + F_3$.   Denote by $r$ the rank of $F$.  	Then
\[m=\mathrm{rank}\, D_2 \leq \mathrm{rank}\, Y_2 +\mathrm{rank}\, F_2 \leq \mathrm{rank}\, Y_2 +r, \]
and similarly
\[
\mathrm{rank}\, Y_3 \leq \mathrm{rank}\, D_3 +\mathrm{rank}\, F_3 \leq r+1.\]
But $\mathrm{rank}\, Y_2 = \mathrm{rank}\, Y_3$, since $Y$ has the common rank property,
and so it follows that
\[m\leq r+r+1,\]
and thus $r\geq \lfloor\frac{m-1}{2}\rfloor$.
Hence,
\[\rho(D, Y)\geq \lfloor\frac{m-1}{2}\rfloor.\]
\end{enumerate}
\end{pf}


\subsection{} \label{sec2.09}
An inspection of the proof of part (b) of the above theorem shows that the finite-dimensionality of the underlying Hilbert space did not really play a role.   In fact, if $\hilb$ is infinite-dimensional, $0 \le j, k < \infty$,  $D \in \bofh$ is normal and $P \in \bofh$ is a projection for which 
\[
\mathrm{rank}\, (I-P) D P = j \mbox{\ \ \ \ \ and \ \ \ \ \ } \mathrm{rank}\, P D (I-P) = k,\] 
then the same argument shows that $\mathrm{rank} (D - Y) \ge \lfloor \frac{|k-j|}{2} \rfloor$ for all operators $Y \in \bofh$ with the (CR) property.


\vskip 2 cm

\section{The infinite-dimensional case}


\subsection{} \label{sec3.01}
Throughout this section, we shall assume that the underlying Hilbert space $\hilb$ is infinite-dimensional and separable.
Our first goal in this section is to extend Theorem~\ref{thm2.05} to this setting.  


\begin{thm} \label{thm3.02}
For all  $0 \le j, k \le \infty$, there exist a normal operator $D \in \bofh$ and an orthogonal projection $P \in \bofh$ for which 
\[
\mathrm{rank}\, (I-P) D P = j \mbox{\ \ \ \ \ and \ \ \ \ \ } \mathrm{rank}\, P D (I-P) = k. \]
\end{thm}


\begin{pf}
By replacing $P$ by $I-P$ if necessary, it becomes clear that there is no loss of generality in assuming that $j \le k$.

\bigskip

\noindent{\textsc{Case One}}:  $j=0$.  

If $k = 0$ as well, we may consider $D = I$, the identity operator, and let $P$ be any non-zero projection.

For $k = 1$, we consider the bilateral shift $U$:  that is, let $\{ e_n\}_{n=1}^\infty$ be an orthonormal basis for $\hilb$, and set $U e_n = e_{n-1}$ for all $n \in \bbZ$.  Let $P_0$ denote the orthogonal projection of $\hilb$ onto $\ol{\mathrm{span}}\, \{ e_n\}_{n \le 0}$.  The condition above is satisfied with  $D := U$, $P = P_0$.

For $2 \le k \le \infty$, we simply consider the tensor product $D:= U \otimes I_k$ of $U$ above with $I_k$, the identity operator acting on a Hilbert space $\cK$ of dimension $k$, and we set  $P = P_0 \otimes I_k$ to obtain the desired rank equalities.  

\bigskip

\noindent{\textsc{Case Two}}:  $1 \le j < \infty$.   

Let $U$ denote the bilateral shift from Case One, and $P_0$ denote the orthogonal projection of $\hilb$ onto $\ol{\mathrm{span}}\, \{ e_n\}_{n \le 0}$.   If $ H := (U + U^*) \otimes I_j$, it is relatively straightforward to verify that with $Q_1 := P_0 \otimes I_j$, we have that 
\[
\mathrm{rank}\, (I-Q_1) H Q_1 = j = \mathrm{rank}\, Q_1 H (I-Q_1). \]
Next, let $R = U \otimes I_{k - j}$ (where $\infty - j := \infty$) and choose a projection $Q_2 = P_0 \otimes I_{k-j}$ as in Case One such that 
\[
\mathrm{rank}\, (I-Q_2) R Q_2 = 0 \mbox{\ \ \ \ \ and \ \ \ \ \ } \mathrm{rank}\, Q_2 R (I-Q_2) = k-j. \]
A routine calculation shows that with $D:= H \oplus R$ and $P := Q_1 \oplus Q_2$, the desired rank equalities are met.

\bigskip

\noindent{\textsc{Case Three}}:  $j = \infty$.   

Since we have reduced the problem to the case where $j \le k$, it  follows that $k = \infty$ as well.  

Consider the selfadjoint operator $\hat{H} = \begin{bmatrix} 1 & 1 \\ 1 & 1 \end{bmatrix} \in \bbM_2(\bbC)$.   Then $H := \hat{H} \otimes I = \begin{bmatrix} I & I \\ I & I \end{bmatrix}$ satisfies the condition relative to the projection $P = I \oplus 0$.

\end{pf}



The case where $j = 1$ and $k = \infty$ in the above Theorem is only one possible infinite-dimensional analogue of Theorem~\ref{thm2.03}.   Alternatively, we may view that Theorem as requiring that $D_2$  be invertible.   Interestingly enough, this is no longer possible in the infinite-dimensional setting.  In fact, a stronger (negative) result holds.

\smallskip

\begin{prop} \label{prop3.03}
There does not exist a normal operator
\[
	D=\begin{bmatrix} D_1 & D_2 \\ D_3 & D_4 \end{bmatrix} \]
in $\cB(\hilb\oplus \hilb)$ such that $D_2$ is invertible and $D_3$ is compact.
\end{prop}


\begin{pf}
We argue by contradiction.   If such normal operator $D$ were to exist, it would follow that
\[D_2 D_2^*=(D_1^* D_1 -D_1  D_1^*)+D_3^* D_3.\]
Since $D_2$ is invertible, $D_2 D_2^*$ is positive and invertible, and thus $0$ is not in the essential numerical range of $D_2 D_2^*$.
On the other hand, by a result of the second author~\cite[Theorem~8]{Rad1966}, and keeping in mind that $D_3$ is compact, $0$ is indeed in the essential numerical range of $(D_1^* D_1 -D_1  D_1^*)+D_3^* D_3$, a contradiction.
\end{pf}


\subsection{} \label{sec3.04}
When $1 \le m < \infty$, it is clear that an operator $D_2 \in \bbM_m(\bbC)$ is invertible if and only if $D_2$ is a \textbf{quasiaffinity}; i.e. it is injective and has dense range.   Moreover, in the infinite-dimensional setting, not every normal operator is orthogonally reductive.   Despite this, in light of Proposition~\ref{prop3.03}, the next example is somewhat surprising.   


\begin{thm} \label{thm3.05}
There exists a normal operator
\[
	D=\begin{bmatrix} D_1 & D_2 \\ 0 & D_4 \end{bmatrix} \]
in $\cB(\hilb\oplus \hilb)$ such that $D_2$ is a quasiaffinity.
\end{thm}


\begin{pf}
Let $A =U+2U^*$ and $B =A^*=U^*+2U$, where $U$ is the bilateral shift operator (i.e. $Ue_n=e_{n-1}$, $n \in \mathbb{Z}$) from Theorem~\ref{thm3.02}.   Then $D := A \oplus B$ is easily seen to be a normal operator.

Let $M \in B(\hilb)$ be a positive contraction, and let $N := (I-M^2)^{1/2}$, so that $M N = N M$ and  $M^2 + N^2 = I$. From this it follows that
\[
P:=
	\begin{bmatrix}
	M^2 & M N \\
	M N & N^2
	\end{bmatrix} \]
is an orthogonal projection in $B(\hilb\oplus \hilb)$.    Arguing as in Theorem~\ref{thm2.03}, we see that $\begin{bmatrix} M \\ N \end{bmatrix}$ and $\begin{bmatrix} N \\ - M \end{bmatrix}$ are both isometries from $\hilb$ into $\hilb\oplus\hilb$, and that it suffices to find $M$ and $N$ as above such that 
\[
(N A M - M B N) = \begin{bmatrix} N & -M \end{bmatrix} \ \begin{bmatrix} A & 0 \\  0 & B \end{bmatrix} \ \begin{bmatrix} M \\ N \end{bmatrix} =0, \]
while
\[
(M A N - N B M) =\begin{bmatrix} M & N \end{bmatrix} \ \begin{bmatrix} A & 0 \\  0 & B \end{bmatrix} \ \begin{bmatrix} N \\ -M \end{bmatrix}\textup{is injective and has dense range}.\]

\smallskip

We shall choose $M$ (and thus $N$) to be diagonal operators relative to the orthonormal basis $\{ e_n\}_{n \in \bbZ}$, $M = \mathrm{diag} (\alpha_n)_{n \in \bbZ}$, where $\alpha_n := \dss \frac{1}{\sqrt{1 + 4^{-n}}}$ for each $n \in \bbZ$.  The condition that $N = (I-M^2)^{1/2}$ implies that $N = \mathrm{diag} (\beta_n)_{n \in \bbZ}$, where $\beta_n = \dss \frac{2^{-n}}{\sqrt{1+4^{-n}}}$ for all $n \in \bbZ$.

It is easy to see that $M$ and $N$ are commutative, positive contractions and $M^2 + N^2 = I$ by construction.

Next, 
\begin{align*}
N A M e_n 
	&= N A (\alpha_n e_n) \\
	&= \alpha_n N (e_{n-1} + 2 e_{n+1}) \\
	&= \alpha_n (\beta_{n-1} e_{n-1} + 2 \beta_{n+1} e_{n+1}), 
\end{align*}
while
\begin{align*}
M A^* N e_n 
	&= M A^* (\beta_n e_n) \\
	&= \beta_n M (e_{n+1} + 2 e_{n-1}) \\
	&= \beta_n (\alpha_{n+1} e_{n+1} + 2 \alpha_{n-1} e_{n-1}).
\end{align*}
But 
\[
\alpha_n \beta_{n-1}  =  \frac{1}{\sqrt{1 + 4^{-n}}} \ \frac{2^{-(n-1)}}{\sqrt{1+4^{-(n-1)}}} =  \frac{2}{\sqrt{1 + 4^{-(n-1)}}} \ \frac{2^{-n}}{\sqrt{1+4^{-n}}} = 2 \alpha_{n-1} \beta_n, \]	
and similarly 
\[
2 \alpha_n \beta_{n+1} = \frac{2}{\sqrt{1+4^{-n}}} \frac{2^{-(n+1)}}{\sqrt{1+4^{-(n+1)}}} = \frac{1}{\sqrt{1+4^{-(n+1)}}} \frac{2^{-n}}{\sqrt{1+4^{-n}}} = \alpha_{n+1} \beta_n. \]
Since this holds for all $n \in \bbZ$,  $N A M - M A^* N = 0$, as claimed.	

\smallskip

As for the second equation we must verify, observe that
\[(M A N - N A^* M)^*=NA^*M-MAN=-(M A N - N A^* M).\]
As such, we need only show that $M A N - N A^* M$ is injective, since then $(M A N - N A^* M)^*$ is also injective and thus both are injective and have dense range.

Again, we compute, for each $n \in \bbZ$, 
\begin{align*}
(M A N - N A^* M) e_n
	&= M A N e_n - N A^* M e_n \\
	&= M A (\beta_n e_n) - N A^* (\alpha_n e_n) \\
	&= \beta_n M (e_{n-1} + 2 e_{n+1}) - \alpha_n N (e_{n+1} + 2 e_{n-1}) \\
	&= \beta_n (\alpha_{n-1} e_{n-1} + 2 \alpha_{n+1} e_{n+1}) - \alpha_n (\beta_{n+1} e_{n+1} + 2 \beta_{n-1} e_{n-1}) \\
	&= (\alpha_{n-1} \beta_n - 2 \alpha_n \beta_{n-1}) e_{n-1} + (2 \alpha_{n+1} \beta_n - \alpha_n \beta_{n+1}) e_{n+1}.
\end{align*}
Suppose that $x = \sum_{n \in \bbZ} x_n e_n \in \ker\, (MAN - N A^* M)$.  	Then 
\begin{align*}
0
	&= (MAN - N A^* M)  \sum_{n \in \bbZ} x_n e_n \\
	&= \sum_{n \in \bbZ}  x_n \left((\alpha_{n-1} \beta_n - 2 \alpha_n \beta_{n-1}) e_{n-1} + (2 \alpha_{n+1} \beta_n - \alpha_n \beta_{n+1}) e_{n+1}\right)
\end{align*}	
By equating coefficients, we see that for all $p \in \bbZ$, 
\[
x_{p+1} (\alpha_p \beta_{p+1} - 2 \alpha_{p+1} \beta_p) + x_{p-1} ( 2 \alpha_p \beta_{p-1} - \alpha_{p-1} \beta_p) = 0,\]
or equivalently,
\[
x_{p+1} = - \frac{2 \alpha_p \beta_{p-1} - \alpha_{p-1} \beta_p}{\alpha_p \beta_{p+1} - 2 \alpha_{p+1} \beta_p}  x_{p-1} \mbox{\ \ \ \ \ for all } p \in \bbZ. \]
But a routine calculation shows that 
\[
\frac{2 \alpha_p \beta_{p-1} - \alpha_{p-1} \beta_p}{\alpha_p \beta_{p+1} - 2 \alpha_{p+1} \beta_p} = -2 \frac{\sqrt{1+ 4^{-(p+1)}}}{\sqrt{1+4^{-(p-1)}}},\]
and so the condition that $\norm x \norm^2 = \sum_{p \in \bbZ} |x_p|^2 < \infty$ clearly implies that 
\[
x_p = 0 \mbox{\ \ \ for all } p \in \bbZ. \]
Thus $\ker\, (M A N - N A^* M) = 0 = \ker\, (M A N - N A^* M)^*$, as required to complete the proof.
\end{pf}


Using a slightly more subtle ``direct sum" device than in Case Two of Theorem~\ref{thm3.02}, we obtain:

\smallskip

\begin{cor} \label{cor3.06}
If $1 \le j$ is any positive integer,  then there exists a normal operator $D \in \cB(\hilb\oplus\hilb)$ and a projection $P\in \cB(\hilb\oplus\hilb)$ of infinite rank and nullity such that 
\[
\mathrm{rank}\, (I-P) D P = j \]
and $P D (I-P)$ is a quasiaffinity.
\end{cor}


\begin{pf}
By Theorem~\ref{thm3.05}, we can find a normal operator $N = \begin{bmatrix} N_1 & N_2 \\ 0 & N_4 \end{bmatrix} \in \cB(\hilb \oplus \hilb)$ such that $N_2$ is a quasiaffinity.   Let $Q = \begin{bmatrix} I_j & I_j \\ I_j & I_j \end{bmatrix} \in \bbM_{2 j}(\bbC)$, so that $Q$ is ($2$ times) a projection of rank $j$.  Then $D := N \oplus Q$ is clearly normal, and it is unitarily equivalent to 
\[
\begin{bmatrix} 
	I_j & & & I_j \\ & N_1 & N_2 & \\ & 0 & N_4 & \\ I_j & & & I_j 
\end{bmatrix}. \]
Set $D_1 = \begin{bmatrix} I_j & 0 \\ 0 & N_1\end{bmatrix}$, $D_2 = \begin{bmatrix} 0 & I_j \\ N_2 & 0 \end{bmatrix}$, $D_3 = \begin{bmatrix} 0 & 0 \\ I_j & 0 \end{bmatrix}$ and $D_4 = \begin{bmatrix} N_4 & 0 \\ 0 & I_j \end{bmatrix}$.   

Clearly $\mathrm{rank}\, D_3 = j$ and $D_2$ is a quasiaffinity.	
\end{pf}


\subsection{} \label{sec3.07}
In Theorem~\ref{thm2.08}, we saw that if $D \in \bbM_{2m}(\bbC)$ is a normal matrix, and if $P \in \bbM_{2m}(\bbC)$ is a projection of rank $m$ such that $\mathrm{rank}\, (I-P)D P = 1$ and $\mathrm{rank}\, P D (I-P) = m$, then $D$ is necessarily cyclic.   It is reasonable to ask, therefore, whether an analogue of this might hold in the infinite-dimensional setting.   In general, the answer is no.


\begin{cor} \label{cor3.08}
For any integer $j \ge 0$, there exists a non-cyclic normal operator $D \in \cB(\hilb)$ and an orthogonal projection $P\in \cB(\hilb)$ of infinite rank and nullity such that 
\[
\mathrm{rank}\, (I-P) D P = j \] 
and $P D (I-P)$ is a quasiaffinity.
\end{cor}


\begin{pf}
By Theorem~\ref{thm3.05}, we can choose a normal operator $N \in \bofh$
\[
	N=\begin{bmatrix} N_1 &N_2\\0 & N_4 \end{bmatrix}, \]
where $N_2$ is a quasiaffinity, and by Corollary~\ref{cor3.06} (or by Theorem~\ref{thm3.05} once again if $j=0$), we may choose a normal operator $M \in \bofh$ such that 
\[
	M=\begin{bmatrix} M_1 & M_2 \\ M_3 & M_4 \end{bmatrix}, \]
where $\mathrm{rank}\, M_2 =j$ and $M_2$ is a quasiaffinity.

Define
\[
	D=\begin{bmatrix} N_1&&&&&&&N_2\\&N_{1}&&&&N_{2}&\\ &&&M_1 & M_2&&& \\ &&&M_3 & M_4&&&\\&0&&&&N_{4}& \\ 0&&&&&&&N_4 \end{bmatrix}.\]

Letting $P = I \oplus I \oplus I \oplus 0 \oplus 0 \oplus 0$, we see that $\mathrm{rank}\, (I-P) D P = \mathrm{rank}\, M_3 = j$ and $P D (I-P)$ is a quasiaffinity.   Moreover, $D$ is unitarily equivalent to $N \oplus N \oplus M$, and thus is not cyclic.
\end{pf}



\vskip 1 cm

\section{Compact normal operators}


\subsection{} \label{sec4.01}
Let $D \in \bofh$ (where $\hilb$ is either finite- or infinite-dimensional) be a normal operator, and let $P \in \bofh$ be a non-trivial projection.   Write 
\[
D = \begin{bmatrix} D_1 & D_2 \\ D_3 & D_4 \end{bmatrix} \]
relative to the decomposition $\hilb = \mathrm{ran} \, P \oplus \mathrm{ran}\, (I-P)$.

The fact that in the infinite-dimensional setting we can find $D$ and $P$ as above such that $D_3 = 0 \ne D_2$, whereas no such $D$ and $P$ exist when $\mathrm{dim}\, \hilb < \infty$ is the statement that not every normal operator acting on an infinite-dimensional Hilbert space is orthogonally reductive, whereas every normal matrix is.

In~\cite{APTT2009}, the concept of an \emph{almost-invariant} subspace for bounded linear operators $T$ acting on infinite-dimensional Banach spaces was introduced.   Given a Banach space $\fX$ and an infinite-dimensional (closed) subspace $\fM$ of $\fX$ such that $\fX/\fM$ is again infinite-dimensional ($\fM$ is then called a \textbf{half-space} of $\fM$), we say that $\fM$ is \textbf{almost-invariant} for $T$ if there exists a finite-dimensional subspace $\fF$ of $\fX$ such that $T \fM \subseteq \fM + \fF$.   The minimal dimension of such a space $\fF$ is referred to as the \textbf{defect} of $T$ relative to $\fM$.    In~\cite{PT2013} and~\cite{Tca2019}, it was shown that every operator $T$ acting on an infinite-dimensional Banach space admits an almost-invariant half-space of defect at most 1.  This is a truly remarkable result.

As a possible generalisation of the notion of reductivity for Hilbert space operators, we propose the following definition.


\begin{defn} \label{defn4.02}
An operator $T \in \bofh$ is said to be \textbf{almost reductive} if for every projection $P \in \bofh$, the condition that $\mathrm{rank}\, (I-P) T P < \infty$ implies that $\mathrm{rank}\, P T (I-P) < \infty$.   
\end{defn}


\subsection{} \label{sec4.03}
It is clear that every invariant-half space is automatically almost-invariant for $T$.   If the notion of ``almost-reductivity" is to make sense, one should expect that every \emph{orthogonally reductive} operator should be ``almost reductive".   

The relevance of this to the problem we have been examining is as follows:   if $K \in \bofh$ is a \emph{compact}, normal operator, then it is well-known~\cite{Wer1952} that $K$ is orthogonally reductive.   This leads to the following question.
 

\begin{ques} \label{ques4.04}
Is every compact, normal operator $K$ almost reductive?  (More generally, is every \emph{reductive} normal operator $D \in \bofh$ almost reductive?)
\end{ques}

\smallskip

Phrased another way, does there exist a compact, normal operator $K$ and a projection $P$ (necessarily of infinite rank and nullity) such that 
\[
\mathrm{rank}\, (I-P) K P < \infty \mbox{\ \ \ \ \ and \ \ \ \ \ } \mathrm{rank}\, P K (I-P) = \infty?\]
The normal operators $D$ constructed in Theorem~\ref{thm3.02} and Theorem~\ref{thm3.05} for which $\mathrm{rank}\, (I-P) D P < \infty$ and $\mathrm{rank}\, P D (I-P) = \infty$ were definitely not compact, and nor were they reductive.

\bigskip

So far, we have been unable to resolve this question.   Indeed, we propose the following (potentially simpler) question:


\begin{ques} \label{ques4.05}
Do there exist a compact, normal operator $K \in \bofh$ and a projection $P \in \bofh$ such that $\mathrm{rank}\, (I-P) K P < \infty$ and $P K (I-P)$ is a quasiaffinity?
\end{ques}
  

While we do not have an answer to this question, nevertheless, there are some things that we can say about its structure, should such an operator $K$ exist.  First we recall a result of Fan and Fong which we shall require.

\begin{thm}~\cite[Theorem 1]{FanFon80} \label{thm4.06}
Let $H$ be a compact, hermitian operator. Then the following are equivalent:
\begin{enumerate}
	\item[(a)] 
	$H=[A^*, A]$ for some compact operator $A$.
	\item[(b)] 
	There exists an orthonormal basis $\{e_n\}_{n\in \mathbb{N}}$ such that $\langle H e_n, e_n \rangle=0$ for all  $n\in \mathbb{N}$.
\end{enumerate}
\end{thm}


\bigskip

Recall that a compact operator $K \in \bofh$ is said to be a \textbf{Hilbert-Schmidt} operator if there exists an orthonormal basis $\{ e_n\}_{n=1}^\infty$ for $\hilb$ such that
\[
\norm K \norm_2 := \left( \mathrm{tr} (K^* K)\right)^{1/2} = \left( \sum_{n=1}^\infty \langle K^* K e_n, e_n\rangle \right) ^{1/2} < \infty.\]
(Equivalently, this holds for \emph{all} orthonormal bases $\{ e_n\}_{n=1}^\infty.)$  When this is the case, the map $K \mapsto \norm K \norm_2$ defines a norm on the set $\cC_2(\hilb)$ of all Hilbert-Schmidt operators on $\hilb$.   (Although this is not the original definition of $\cC_2(\hilb)$, it is equivalent to it.)

\smallskip

\begin{cor} \label{cor4.07}
Let
\[
	K =\begin{bmatrix} K_1 & K_2 \\ K_3 & K_4 \end{bmatrix} \]
be a compact, normal operator in $\cB(\hilb\oplus \hilb)$.
Then $K_2 \in C_2(\hilb)$ if and only if $K_3 \in C_2(\hilb)$, in which case $\| K_2 \|_2=\| K_3 \|_2$.

In particular, therefore, if $K_3$ is a finite-rank operator, then $K_2$ must be a Hilbert-Schmidt operator.
\end{cor}


\begin{pf}
Since $K$ is normal, it follows that  $K_1^* K_1 +K_3^* K_3 =K_1 K_1^*+K_2 K_2^*$, and thus $[K_1^*, K_1] = K_2 K_2^*-K_3^*K_3$.
Now $K_1$ is compact, and so by the above theorem, there exists an orthonormal basis $\{e_n\}_{n\in \mathbb{N}}$ such that $\langle (K_2 K_2^* - K_3^*K_3)e_n, e_n \rangle =0$ for all $n \in \mathbb{N}$. 

Suppose that $K_3\in C_2(\hilb)$.  Then
\[
\|K_3 \|_2^2=\textup{tr}(K_3^*K_3)=\sum_{n =1}^\infty \langle K_3^*K_3 e_n, e_n \rangle<\infty.\]
Therefore,
\[\sum_{j\in \mathbb{N}}\langle K_2 K_2^* e_n, e_n \rangle =\sum_{n=1}^\infty \langle K_3^* K_3 e_n, e_n \rangle<\infty,\]
proving that $K_2 \in C_2(\hilb)$, and
\[\|K_2 \|_2 =\| K_3 \|_2.\]

The last statement is obvious.
\end{pf}


\subsection{} \label{sec4.08}
The proof of Theorem~\ref{thm2.08} yields a very specific structure result for normal matrices $D \in \bbM_{2m}(\bbC)$ for which there exists an orthogonal projection $P$ satisfying $\mathrm{rank}\, (I-P) D P = 1$ and $\mathrm{rank}\, P D (I-P) = m$.    Since orthogonal reductivity and normality of matrices coincide,  Proposition~\ref{prop4.10} below can be seen as an extension of that structure result to the infinite-dimensional setting.


\begin{defn} \label{defn4.09}
By a \textbf{simple bilateral chain of subspaces} of a Hilbert space $\hilb$ we mean a sequence of closed subspace $\{\cM_j\}_{j=-\infty}^{\infty}$ with
\[\cdots \subset \cM_{-2} \subset \cM_{-1}\subset \cM_0\subset \cM_1\subset \cM_2 \subset \cdots, \]
where $\dim(\cM_{j+1}\ominus \cM_j)=1$ for all $j\in \mathbb{Z}$.
We say an operator $T\in B(\hilb)$ \textbf{shifts forward} a simple bilateral chain $\{\cM_j\}_{j=-\infty}^{\infty}$ if
\[T\cM_j\subset \cM_{j+1},~~\forall j\in \mathbb{Z}.\]
\end{defn}


\begin{prop} \label{prop4.10}
Let $T$ be an orthogonally reductive operator on $\hilb$ and assume that
relative to a decomposition $\hilb=\hilb_1\oplus \hilb_2$, it has the representation
\[T =\begin{bmatrix} A & L \\ F & B \end{bmatrix}, \]
where $F$ has rank one and $L$ has infinite rank. Then
$T$ has an infinite-dimensional invariant subspace $\hilb_0$
such that the restriction $T_0$ of $T$ to $\hilb_0$ shifts forward a simple bilateral chain $\{\cM_j\}_{j=-\infty}^{\infty}$ of subspaces.
\end{prop}


\begin{pf}
Assume with no loss that $T$ is invertible and let $\cM_0=\hilb_1$. We will define subspace $\cM_j$ inductively:  we set

\[\cM_{j+1}=\cM_j+T\cM_j,~~\forall j\geq 0\]
and
\[\cM_{j-1}=\cM_j\cap T^{-1}\cM_j,~~\forall j\leq 0.\]
Then
\[\cdots\subset \cM_{-1}\subset \cM_0\subset \cM_1\subset \cdots,\]
and
\[T\cM_j\subset \cM_{j+1},~~\forall j\in \mathbb{Z}.\]
The assumption that $F$ has rank one implies that $\cM_1\ominus \cM_0$ has dimension one.  It follows inductively that
the dimension of $\cM_{j+1}\ominus \cM_j$ is at most one for all $j\in \mathbb{Z}$.
We shall show that this difference in dimensions is exactly one for all $j\in \mathbb{Z}$.

Suppose not. First assume $j>1$.
If $\cM_{j+1}=\cM_j$, then $\cM_j$ is invariant under $T$ and thus reducing. This means that
\[P_jT(I-P_j)=0, \]
with $P_j$ denoting the orthogonal projection onto $\cM_j$.
In particular, then
\[L(I-P_j)=P_0L(I-P_j)=0.\]
But this implies that the rank of $L$ is at most $j$, which is a contradiction.

The proof for $j<1$ is similar. In summary, we conclude that
$\{\cM_j\}_{j=-\infty}^{\infty}$ is a proper bilateral chain of subspaces.

Now $\bigcap_{j=-\infty}^{\infty}\cM_j$ and $\bigvee_{j=-\infty}^{\infty}\cM_j$ are both invariant,
and hence reducing. Let
\[\hilb_0=(\bigvee_{j=-\infty}^{\infty}\cM_j)\ominus (\bigcap_{j=-\infty}^{\infty}\cM_j),\]
and note that if we define
\[\cM_j'=\cM_j\ominus (\bigcap_{k=-\infty}^{\infty}\cM_k),\]
and $T_0 := T|_{\hilb_0}$,
then
 $\{\cM_j'\}_{j=-\infty}^{\infty}$ is the desired bilateral chain in $\hilb_0$  which $T_0$ shifts forward.
\end{pf}


For compact normal operators, we can obtain a stronger result.

\begin{cor} \label{cor4.11}
If $K$ is a compact normal operator on $\hilb=\hilb_1\oplus \hilb_2$ of the form
\[K =\begin{bmatrix} A & L \\ F & B \end{bmatrix}, \]
where $F$ has rank one and $L$ is a quasiaffinity, then
$K$ shifts forward a simple bilateral chain $\{\cM_j\}_{j=-\infty}^{\infty}$ of subspaces.  $($Here it is understood that $\mathrm{dim}\, \hilb_1 = \infty = \mathrm{dim}\, \hilb_2$.$)$
\end{cor}


\begin{pf}
It is well-known that compact normal operators are orthogonally reductive \cite{Wer1952}. Thus we must only show that the subspace $\hilb_0$ of the
proposition above coincides with $\hilb$. In other words,
\[\bigcap_{j=-\infty}^{\infty}\cM_j=0,~~\bigvee_{j=-\infty}^{\infty}\cM_j=\hilb.\]

Set
\[\mathcal{N}_1=\bigcap_{j=-\infty}^{\infty}\mathcal{M}_j,~~\mathcal{N}_2=\mathcal{M}_0\ominus \mathcal{N}_1,~~
\mathcal{N}_3=(\bigvee_{j=-\infty}^{\infty}\mathcal{M}_j)\ominus \mathcal{M}_0,~~\mbox{ and }  \mathcal{N}_4=(\bigvee_{j=-\infty}^{\infty}\mathcal{M}_j)^{\perp}.\]
As $\mathcal{N}_1, \bigoplus_{1\leq i\leq 3}\mathcal{N}_i$ are both invariant and therefore reducing for $K$, with respect to the decomposition of $\hilb=\hilb_1\oplus \hilb_2=(\mathcal{N}_1\oplus\mathcal{N}_2)\oplus(\mathcal{N}_3\oplus\mathcal{N}_4)$, we may write
\[
K=\begin{matrix}\begin{bmatrix}
A&L\\
F&B
\end{bmatrix}&\begin{matrix}
  \hilb_1\\
  \hilb_2\\
\end{matrix}
\end{matrix}=\begin{matrix}\begin{bmatrix}
 A_1&0&0&0\\
 0&A_2&L'&0\\
 0&F'&B_3&0\\
 0&0&0&B_4\\
\end{bmatrix}&\begin{matrix}
  \mathcal{N}_1\\
  \mathcal{N}_2\\
  \mathcal{N}_3\\
  \mathcal{N}_4\\
\end{matrix}
\end{matrix}.
\]
Since
\[L=\begin{bmatrix}
0&0\\
L'&0
\end{bmatrix}\]
is a quasiaffinity, it follows that $\mathcal{N}_1=0$, and similarily $\mathcal{N}_4=0$.
In other words,
\[\bigcap_{j=-\infty}^{\infty}\mathcal{M}_j=\mathcal{N}_1=0,~~\bigvee_{j=-\infty}^{\infty}\mathcal{M}_j=\mathcal{N}_4^{\perp}=\hilb.\]

\end{pf}







\end{document}